\documentclass{article}
\usepackage{amsfonts}
\usepackage{amsmath}

\newtheorem{theorem}{Theorem}

\newtheorem{definition}[theorem]{Definition}
\newtheorem{example}[theorem]{Example}

\input{tcilatex}
\begin{document}

\title{\textbf{Inverse Fourier Transform for Bi-Complex Variables}}
\author{A. BANERJEE$^{1}$, S. K. DATTA$^{2}$ and MD. A.\ HOQUE$^{3}$ \\
$^{1}$Department of Mathematics, Krishnath College, Berhampore, \\
\ Murshidabad 742101, India, E-mail: abhijit.banerjee.81@gmail.com\\
$^{2}$Department of Mathematics, University of Kalyani, Kalyani, Nadia,\\
\ PIN-741235,India, E-mail: sanjib\_kr\_datta@yahoo.co.in \ \ \\
$^{3}$Sreegopal Banerjee College,Bagati, Mogra, Hooghly, \\
\ 712148, India,E-mail: mhoque3@gmail.com}
\date{}
\maketitle

\begin{abstract}
In this paper we examine the existence of bicomplexified inverse Fourier
transform \ as an extension of it's complexified inverse version within the
region of convergence of bicomplex Fourier transform. In this paper we use
the idempotent representation of bicomplex-valued functions as projections
on the auxiliary complex spaces of the components of bicomplex numbers along
two orthogonal,idempotent hyperbolic directions.

\end{abstract}

\section{\protect\bigskip Introduction}

In 1892, in search for special algebras, Corrado Segre $\left\{ \text{cf. 
\cite{Segre}}\right\} $ published a paper in which he treated an infinite
family of algebras whose elements are commutative generalization of complex
numbers called bicomplex numbers, tricomplex numbers,.....etc. Segre $%
\left\{ \text{cf. \cite{Segre}}\right\} $ defined a bicomplex number $\xi $
as follows:

\begin{equation*}
\xi=x_{1}+i_{1}x_{2}+i_{2}x_{3}+i_{1}i_{2}x_{4},
\end{equation*}

where $x_{1},$ $x_{2},x_{3},x_{4}$ are real numbers, $i_{1}^{2}=i_{2}^{2}=-1$
and $i_{1}i_{2}=i_{2}i_{1}.$ The set of bicomplex numbers, complex numbers
and real numbers are respectively denoted by $%
\mathbb{C}
_{2},$ $%
\mathbb{C}
_{1}$ and $%
\mathbb{C}
_{0}$ . Thus

\begin{equation*}
\mathbb{C}
_{2}=\{\xi:\xi=x_{1}+i_{1}x_{2}+i_{2}x_{3}+i_{1}i_{2}x_{4},\text{ }%
x_{1},x_{2},x_{3},x_{4}\in%
\mathbb{C}
_{0}\}
\end{equation*}

\begin{equation*}
\text{i.e., }%
\mathbb{C}
_{2}=\{\xi=z_{1}+i_{2}z_{2}:\text{ }%
z_{1}(=x_{1}+i_{1}x_{2}),z_{2}(=x_{3}+i_{1}x_{4})\in%
\mathbb{C}
_{1}\}.
\end{equation*}

There are \ two non trivial elements $e_{1}=\frac{1+i_{1}i_{2}}{2}$ and $%
e_{2}=\frac{1-i_{1}i_{2}}{2}$ \ in $%
\mathbb{C}
_{2}$\ with the properties $\ e_{1}^{2}=e_{1},e_{2}^{2}=e_{2},e_{1}\cdot
e_{2}=e_{2}\cdot e_{1}=0$ and $\ e_{1}+e_{2}=1$ which means that $\ e_{1}$
and $e_{2}$ are idempotents alternatively called orthogonal idempotents. By
the help of the idempotent elements $\ e_{1}$ and $e_{2},$ any bicomplex
number%
\begin{equation*}
\xi
=a_{0}+i_{1}a_{1}+i_{2}a_{2}+i_{1}i_{2}a_{3}=(a_{0}+i_{1}a_{1})+i_{2}(a_{2}+i_{1}a_{3})=z_{1}+i_{2}z_{2}
\end{equation*}%
\ where \ $a_{0,}$ $a_{1},a_{2},a_{3}$ $\in 
\mathbb{C}
_{0},$ 
\begin{equation*}
z_{1}(=a_{0}+i_{1}a_{1})\text{ and }z_{2}(=a_{2}+i_{1}a_{3})\text{ }\in 
\mathbb{C}
_{1}
\end{equation*}%
can be expressed as $\ $%
\begin{equation*}
\xi =z_{1}+i_{2}z_{2}=\xi _{1}e_{1}+\xi _{2}e_{2}
\end{equation*}%
where $\xi _{1}(=z_{1}-i_{1}z_{2})\in 
\mathbb{C}
_{1}$ and $\xi _{2}(=z_{1}+i_{1}z_{2})$ $\in 
\mathbb{C}
_{1}.$

\section{Fourier Transform}

Let $f(t)$ be a real valued continuous function in \ $(-\infty ,\infty )$
which satisfies the estimates%
\begin{align}
|f(t)|& \leq C_{1}\exp (-\alpha t),t\geq 0,\alpha >0  \notag \\
\text{and }|f(t)|& \leq C_{2}\exp (-\beta t),t\leq 0,\beta >0.  \label{6.4.5}
\end{align}%
Then the bicomplex Fourier transform $\left\{ \text{cf. \cite{Banerjee}}%
\right\} $ of $f(t)$ can be defined as

\begin{equation*}
\widehat{f}(\omega )=\ \tciFourier \{f(t)\}=\tint\limits_{-\infty }^{\infty
}\exp (i_{1}\omega t)f(t)dt,\omega \in 
\mathbb{C}
_{2}.
\end{equation*}%
The Fourier transform $\widehat{f}(\omega )$\ exists and holomorphic for all 
$\omega \in $ $\Omega $\ where

\begin{equation*}
\Omega =\{\omega =a_{0}+i_{1}a_{1}+i_{2}a_{2}+i_{1}i_{2}a_{3}\in 
\mathbb{C}
_{2}:-\infty <a_{0},a_{3}<\infty ,
\end{equation*}

\begin{equation*}
-\alpha +|a_{2}|<a_{1}<\beta -|a_{2}|\text{ }and\text{ }0\leq |a_{2}|<\frac{%
\alpha +\beta }{2}\}
\end{equation*}%
is the region of absolute convergence of $\widehat{f}(\omega ).$

\subsection{Complex version of Fourier inverse transform.}

We start with the complex version of Fourier inverse transform and in this
connection we consider a continuous function $f(t)$ for $-\infty<t<\infty$
satisfying \ the estimates ($\ref{6.4.5})$ possessing the Fourier transform $%
\widehat{f_{1}\text{ }}$in complex variable $\omega_{1}$=$x_{1}+i_{1}x_{2}$
i.e.,

\begin{align*}
\widehat{f_{1}}(\omega_{1}) & =\tint
\limits_{-\infty}^{\infty}\exp(i_{1}\omega_{1}t)f(t)dt \\
& =\tint
\limits_{-\infty}^{\infty}\exp(i_{1}x_{1}t)\{\exp(-x_{2}t)f(t)\}dt=%
\phi(x_{1},x_{2}).
\end{align*}

In fact, one may identify $\phi(x_{1},x_{2})$ as the Fourier transform of $%
g(t)=\exp(-x_{2}t)f(t)$ in usual complex exponential form $\left\{ \text{cf. 
\cite{Boc} \& }\cite{Kai}\right\} $.

\ \ Towards this end, we assume that $f(t)$ is continuous and $f^{\prime}(t)$
is piecewise continuous on the whole real line. Then $\widehat{f_{1}}%
(\omega_{1})$\ converges absolutely for $-\alpha<$ $x_{2}<\beta$ and

\begin{equation*}
|\widehat{f_{1}}(\omega_{1})\ |<\infty
\end{equation*}

which implies that

\begin{equation*}
\tint \limits_{-\infty}^{\infty}|\exp(i_{1}\omega_{1}t)f(t)|dt
\end{equation*}

\begin{equation*}
=\tint \limits_{-\infty}^{\infty}|\exp(i_{1}x_{1})g(t)|dt
\end{equation*}

\begin{equation*}
=\tint \limits_{-\infty}^{\infty}|g(t)|dt<\infty.
\end{equation*}

The later condition shows $g(t)$ is absolutely integrable. Then by the
Fourier inverse transform in complex exponential form $\left\{ \text{cf. 
\cite{Boc} \& }\cite{Kai}\right\} ,$

\begin{equation*}
g(t)=\frac{1}{2\pi }\tint\limits_{-\infty }^{\infty }\exp (-i_{1}x_{1}t)\phi
(x_{1},x_{2})dx_{1}
\end{equation*}

which implies that

\begin{equation*}
f(t)=\frac{1}{2\pi}\tint
\limits_{-\infty}^{\infty}\exp(x_{2}t)\exp(-i_{1}x_{1}t)%
\phi(x_{1},x_{2})dx_{1}.
\end{equation*}

Now if we integrate along a horizontal line then $x_{2}$ is constant and so
for complex variable $\omega_{1}$=$x_{1}+i_{1}x_{2}$ (which implies $%
d\omega_{1}=dx_{1}$), the above inversion formula can be extended upto
complex Fourier inverse transform%
\begin{align}
f(t) & =\frac{1}{2\pi}\tint
\limits_{-\infty}^{\infty}\exp\{-i_{1}(x_{1}+i_{1}x_{2})t\}%
\phi(x_{1},x_{2})dx_{1}  \notag \\
& =\frac{1}{2\pi}\tint
\limits_{-\infty+i_{1}x_{2}}^{\infty+i_{1}x_{2}}\exp(-i_{1}\omega_{1}t)%
\widehat{f_{1}}(\omega_{1})d\omega_{1}  \notag \\
& =\frac{1}{2\pi}\lim_{x_{1}\longrightarrow\infty}\tint
\limits_{-x_{1}+i_{1}x_{2}}^{x_{1}+i_{1}x_{2}}\exp(-i_{1}\omega_{1}t)%
\widehat{f_{1}}(\omega_{1})d\omega_{1}.  \label{6.10}
\end{align}
Here the integration is to be performed along a horizontal line in complex $%
\omega_{1}$-plane employing contour integration method.

We first consider the case\textbf{\ }$Im(\omega_{1})=x_{2}\geq0.$We observe
that $\widehat{f_{1}}(\omega_{1})$ is continuous for $x_{2}\geq0$ and in
particular it is holomorphic (and so it has no singularities) for $0$\ $\leq
x_{2}<\beta$ . We now introduce a contour $\Gamma_{R}$\ consisting of the
segment $[-R,R]$ and a semicircle $C_{R}$ of radius $|\omega_{1}|=R>\beta$
with centre at the origin. Then all possible singularities (if exists) of $%
\widehat{f_{1}}(\omega_{1})$ must lie in the region above the horizontal
line $x_{2}=\beta$ . At this stage we now consider the following two cases:

\textbf{CaseI:}We assume that $\widehat{f_{1}}(\omega _{1})$ is holomorphic
in $x_{2}>\beta $ except for having a finite number of poles $\omega
_{1}^{(k)}$ \ for $k=1,2,...n$ therein (See Figure 2 in Appendix). By taking 
$R\rightarrow \infty ,$\ we can guarantee that all these poles lie inside
the contour $\Gamma _{R}$. Since $exp(-i_{1}\omega _{1}t)$ never vanishes
then the status of these poles $\omega _{1}^{(k)}$ of $\widehat{f_{1}}%
(\omega _{1})$ is not affected by multiplication of it with $%
exp(-i_{1}\omega _{1}t)$.Then by Cauchy's residue theorem,%
\begin{align}
& \lim_{R\longrightarrow \infty }\int_{\Gamma _{R}}\exp (-i_{1}\omega _{1}t)%
\widehat{f_{1}}(\omega _{1})d\omega _{1}  \notag \\
& =2\pi i_{1}\sum_{Im(\omega _{1}^{(k)})>0}\func{Re}s\{\exp (-i_{1}\omega
_{1}t)\widehat{f_{1}}(\omega _{1}):\omega _{1}=\omega _{1}^{(k)}\}.
\label{6.11}
\end{align}

Furthermore as $x_{2}\geq 0,$ we can get $|\exp (-i_{1}\omega _{1}t)|\leq 1$
\ for $\omega _{1}\in C_{R}$ only when $t\leq 0$. In particular for $t<0,$

\begin{align*}
M(R) & =\max_{\omega_{1}\in C_{R}}|\widehat{f_{1}}(\omega_{1})|=\max
_{\omega_{1}\in C_{R}}|\tint
\limits_{-\infty}^{0}\exp(i_{1}\omega_{1}t)f(t)dt| \\
& \leq C_{2}\max_{\omega_{1}\in C_{R}}|\tint
\limits_{-\infty}^{0}\exp\{(\beta+i_{1}\omega_{1})t\}dt|=C_{2}\max_{%
\omega_{1}\in C_{R}}|\frac {1}{\beta+i_{1}\omega_{1}}| \\
& \leq C_{2}\max_{\omega_{1}\in C_{R}}\frac{1}{\beta+|i_{1}||\omega_{1}|}
\end{align*}
where we use the estimate $\ref{6.4.5}$. Now for \TEXTsymbol{\vert}$%
\omega_{1}|=R\rightarrow\infty$ , we obtain that $M(R)\rightarrow0$. Thus
the conditions for Jordan's lemma $\left\{ \text{cf. \cite{Sid}}\right\} $
are met and so employing it we get that%
\begin{equation}
\lim_{R\longrightarrow\infty}\int_{C_{R}}\exp(-i_{1}\omega_{1}t)\widehat {%
f_{1}}(\omega_{1})d\omega_{1}=0.  \label{6.12}
\end{equation}
Finally as,

\begin{equation*}
\lim_{R\longrightarrow\infty}\int_{\Gamma_{R}}\exp(-i_{1}\omega_{1}t)%
\widehat{f_{1}}(\omega_{1})d\omega_{1}
\end{equation*}

\begin{equation*}
=\int_{C_{R}}\exp(-i_{1}\omega_{1}t)\widehat{f_{1}}(\omega_{1})d\omega_{1}+%
\tint \limits_{-R+i_{1}x_{2}}^{R+i_{1}x_{2}}\exp(-i_{1}\omega_{1}t)\widehat{%
f_{1}}(\omega_{1})d\omega_{1}
\end{equation*}
then for $R\rightarrow\infty,$ on using ($\ref{6.11})$ and ($\ref{6.12})$ we
obtain that

\begin{equation*}
\tint \limits_{-\infty+i_{1}x_{2}}^{\infty+i_{1}x_{2}}\exp(-i_{1}\omega_{1}t)%
\widehat{f_{1}}(\omega_{1})d\omega_{1}
\end{equation*}

\begin{equation*}
=2\pi i_{1}\sum_{Im(\omega _{1}^{(k)})>0}\func{Re}s\{\exp (-i_{1}\omega
_{1}t)\widehat{f_{1}}(\omega _{1}):\omega _{1}=\omega _{1}^{(k)}\}\text{ for 
}t<0
\end{equation*}

and so

\begin{equation*}
\text{ }f(t)=i_{1}\sum_{Im(\omega_{1}^{(k)})>0}\func{Re}s\{\exp
(-i_{1}\omega_{1}t)\widehat{f_{1}}(\omega_{1}):\omega_{1}=\omega_{1}^{(k)}\}%
\text{ for }t<0.
\end{equation*}

\textbf{Case II: }Suppose $\widehat{f_{1}}(\omega _{1})$ has infinitely many
poles $\omega _{1}^{(k)}$ \ for $k=1,2,...n$ in $x_{2}>\beta $ and we
arrange them in such a way that \TEXTsymbol{\vert}$\omega _{1}^{(1)}$%
\TEXTsymbol{\vert}$\leq |\omega _{1}^{(2)}|\leq |\omega _{1}^{(3)}|.....$%
where \TEXTsymbol{\vert} $\omega _{1}^{(k)}|\rightarrow \infty $ as $%
k\rightarrow \infty $. We then consider a sequence of contours $\Gamma _{k}$
consisting of the segments $[-$ $x_{1}^{(k)}+i_{1}x_{2},$ $%
x_{1}^{(k)}+i_{1}x_{2}]$ and the semicircles $C_{k}$ of radii r$_{k}$ = 
\TEXTsymbol{\vert}$\omega _{1}^{(k)}|>\beta $ enclosing the first k poles $%
\omega _{1}^{(1)}$, $\omega _{1}^{(2)},\omega _{1}^{(3)},.......\omega
_{1}^{(k)}$ (See Figure 3 in Appendix).Then by Cauchy's residue theorem we
get that

\begin{equation*}
2\pi i_{1}\sum_{Im(\omega_{1}^{(k)})>0}\func{Re}\text{s}\{\exp
(-i_{1}\omega_{1}t)\widehat{f_{1}}(\omega_{1}):\omega_{1}=\omega_{1}^{(k)}\}
\end{equation*}

\begin{equation*}
=\int_{\Gamma_{R}}\exp(-i_{1}\omega_{1}t)\widehat{f_{1}}(\omega_{1})d%
\omega_{1}
\end{equation*}

\begin{equation*}
=\int_{C_{R}}\exp(-i_{1}\omega_{1}t)\widehat{f_{1}}(\omega_{1})d\omega_{1}
\end{equation*}

\begin{equation}
+\tint\limits_{-x_{1}^{(k)}+i_{1}x_{2}}^{x_{1}^{(k)}+i_{1}x_{2}}\exp
(-i_{1}\omega _{1}t)\widehat{f_{1}}(\omega _{1})d\omega _{1}.  \label{6.13}
\end{equation}%
Now for t \TEXTsymbol{<} 0, in the procedure similar to Case I, employing
Jordan lemma here also we may deduce that%
\begin{equation*}
\lim_{|\omega _{1}^{(k)}|\longrightarrow \infty }\int_{C_{R}}\exp
(-i_{1}\omega _{1}t)\widehat{f_{1}}(\omega _{1})d\omega _{1}=0.
\end{equation*}

Hence in the limit \ $|\omega _{1}^{(k)}|\longrightarrow \infty $ ( which
implies that $|x_{1}^{(k)}|\longrightarrow \infty )$ , ($\ref{6.13})$ leads
to%
\begin{align*}
& \tint\limits_{-\infty +i_{1}x_{2}}^{\infty +i_{1}x_{2}}\exp (-i_{1}\omega
_{1}t)\widehat{f_{1}}(\omega _{1})d\omega _{1} \\
& =2\pi i_{1}\sum_{Im(\omega _{1}^{(k)})>0}\func{Re}\text{s}\{\exp
(-i_{1}\omega _{1}t)\widehat{f_{1}}(\omega _{1}):\omega _{1}=\omega
_{1}^{(k)}\}\text{ for }t<0
\end{align*}%
and as its consequence%
\begin{equation*}
\text{ }f(t)=i_{1}\sum_{Im(\omega _{1}^{(k)})>0}\func{Re}\text{s}\{\exp
(-i_{1}\omega _{1}t)\widehat{f_{1}}(\omega _{1}):\omega _{1}=\omega
_{1}^{(k)}\}\text{ for }t<0.
\end{equation*}%
Thus for $x_{2}\geq 0$\ , whatever the number of poles is finite or
infinite, from the above two cases we obtain the complex version of Fourier
inverse transform as 
\begin{equation}
\text{ }f(t)=i_{1}\sum_{Im(\omega _{1}^{(k)})>0}\func{Re}\text{s}\{\exp
(-i_{1}\omega _{1}t)\widehat{f_{1}}(\omega _{1}):\omega _{1}=\omega
_{1}^{(k)}\}\text{ for }t<0\text{.}  \label{6.14}
\end{equation}%
We now consider the Case Im\textbf{\ (}$\omega _{1}$\textbf{) = }$x_{2}\leq
0 $. The complex valued function $\widehat{f_{1}}(\omega _{1})$ is
continuous for $x_{2}\leq 0$ and holomorphic in $-\alpha <x_{2}\leq 0$.
Introducing a contour $\Gamma _{R^{\prime }}^{\prime }$ consisting of the
segment $[-R^{\prime },R^{\prime }]$ and a semicircle $C_{R^{\prime
}}^{\prime }$ of radius \TEXTsymbol{\vert}$\omega _{1}|=R^{\prime }>\alpha $
with centre at the origin, we see that all possible singularities (if
exists) of $\widehat{f_{1}}(\omega _{1})$ must lie in the region below the
horizontal line $x_{2}=-\alpha $ . If $\overline{\omega }_{1}^{(k)}$ for $%
k=1,2...$ are the poles in $x_{2}<\alpha $, whatever the number of poles are
finite or not for $R^{\prime }\rightarrow \infty $, in similar to the
previous consideration for $x_{2}\geq 0$ we see that for t \TEXTsymbol{>} 0
the conditions for Jordan lemma are met and so%
\begin{equation}
\text{ }f(t)=-i_{1}\sum_{Im(\omega _{1}^{(k)})<0}\func{Re}\text{s}\{\exp
(-i_{1}\omega _{1}t)\widehat{f_{1}}(\omega _{1}):\omega _{1}=\overline{%
\omega }_{1}^{(k)}\}\text{ for }t>0\text{.}  \label{6.15}
\end{equation}%
We then assign the value of f(t) at t = 0 fulfilling the requirement of
continuity of it in $-\infty <t<\infty $. This completes our procedure in
complex $\omega _{1}$ plane.

Similarly in $\omega _{2}(=y_{1}+i_{1}y_{2})$ plane the complex version of
Fourier inverse transform of $\widehat{f_{2}}(\omega _{2})$ will be 
\begin{equation}
f(t)=\frac{1}{2\pi }\lim_{y_{1}\longrightarrow \infty
}\tint\limits_{-y_{1}+i_{1}y_{2}}^{y_{1}+i_{1}y_{2}}\exp (-i_{1}\omega _{2}t)%
\widehat{f_{2}}(\omega _{2})d\omega _{2}\text{ }  \label{6.16}
\end{equation}%
where the integration is to be performed along the horizontal line in $%
\omega _{2}$ plane. Employing \ the contour integration method, we can
obtain that%
\begin{align}
f(t)& =i_{1}\sum_{Im(\omega _{2}^{(k)})>0}\func{Re}\text{s}\{\exp
(-i_{1}\omega _{2}t)\widehat{f_{2}}(\omega _{2}):\omega _{2}=\omega
_{2}^{(k)}\}\text{ for }t<0  \notag \\
& =-i_{1}\sum_{Im(\omega _{2}^{(k)})<0}\func{Re}\text{s}\{\exp (-i_{1}\omega
_{2}t)\widehat{f_{2}}(\omega _{2}):\omega _{2}=\omega _{2}^{(k)}\}\text{ for 
}t>0  \label{6.17}
\end{align}%
and the value of f(t) at t = 0 can be assigned fulfilling the requirement of
continuity of it in $-\infty <t<\infty $.

\subsection{Bicomplex version of Fourier inverse transform.}

Suppose $\widehat{f}(\omega)$ is the bicomplex Fourier transform of the real
valued continuous function f(t) for $-\infty<t<\infty$ where $\omega
=\omega_{1}e_{1}+\omega_{2}e_{2}$ and $\widehat{f}(\omega)=\widehat{f}%
_{1}(\omega_{1})e_{1}+\widehat{f}_{2}(\omega_{2})e_{2}$ in their idempotent
representations. Here the symbols $\omega_{1},\omega_{2},\widehat{f}_{1}$
and $\widehat{f}_{2}$ have their same representation as defined in
Subsection 6.4.1. Then $\widehat{f}(\omega)$ is holomorphic in%
\begin{align}
\Omega & =\{\omega=(x_{1}+i_{1}x_{2})e_{1}+(y_{1}+i_{1}y_{2})e_{2}\in%
\mathbb{C}
_{2}  \notag \\
& :-\alpha<x_{2},y_{2}<\beta,-\infty<x_{1},y_{1}<\infty\}\text{.}
\label{6.18}
\end{align}
Now using complex inversions $\ref{6.10}$ and $\ref{6.16},$ we obtain that%
\begin{align}
f(t) & =f(t)e_{1}+f(t)e_{2}  \notag \\
& =[\frac{1}{2\pi}\tint \limits_{D_{1}}\exp(-i_{1}\omega_{1}t)\widehat{f}%
_{1}(\omega_{1})d\omega_{1}]e_{1}+[\frac {1}{2\pi}\tint
\limits_{D_{2}}\exp(-i_{1}\omega_{2}t)\widehat{f}_{2}(\omega_{2})d%
\omega_{2}]e_{2}  \notag \\
& =\frac{1}{2\pi}\tint
\limits_{D}\exp\{-i_{1}(\omega_{1}e_{1}+\omega_{2}e_{2})t\}\{\widehat{f}%
_{1}(\omega _{1})e_{1}+\widehat{f}_{2}(\omega_{2})e_{2}\}d(\omega_{1}e_{1}+%
\omega_{2}e_{2})  \notag \\
& =\frac{1}{2\pi}\tint \limits_{D}\exp\{-i_{1}(\omega t)\widehat{f}%
(\omega)d\omega\text{ }  \label{6.19}
\end{align}
where

\begin{equation*}
D_{1}=\{\omega=x_{1}+i_{1}x_{2}\in%
\mathbb{C}
(i_{1}):-\infty<x_{1}<\infty,-\alpha<x_{2}<\beta\},
\end{equation*}

\begin{equation*}
D_{2}=\{\omega=y_{1}+i_{1}y_{2}\in%
\mathbb{C}
(i_{1}):-\infty<y_{1}<\infty,-\alpha<y_{2}<\beta\}
\end{equation*}
and D be such that $D_{1}=P_{1}(D),D_{2}=P_{2}(D)$. The integration in $%
D_{1} $ and $D_{2}$ are to be performed along the lines parallel to $x_{1}$
-axis in $\omega_{1}$ plane and $y_{1}$-axis in $\omega_{2}$ plane
respectively inside the respective strips $-\alpha<x_{2}<\beta$ and $%
-\alpha<y_{2}<\beta$. As a result,%
\begin{equation}
D=\{\omega\in%
\mathbb{C}
_{2}:\omega=\omega_{1}e_{1}+%
\omega_{2}e_{2}=(x_{1}+i_{1}x_{2})e_{1}+(y_{1}+i_{1}y_{2})e_{2}\}\text{ \ \
\ }  \label{6.20}
\end{equation}
where $-\infty<x_{1},y_{1}<\infty,-\alpha<x_{2},y_{2}<\beta.$ In
four-component form D can be alternatively expressed as 
\begin{align*}
D & =\{\omega\in%
\mathbb{C}
_{2}:\frac{x_{1}+y_{1}}{2}+i_{1}\frac{x_{2}+y_{2}}{2}+i_{2}\frac{y_{2}-x_{2}%
}{2}+i_{1}i_{2}\frac{x_{1}-y_{1}}{2}, \\
-\infty & <x_{1},y_{1}<\infty,-\alpha<x_{2},y_{2}<\beta\}.
\end{align*}

Conversely, if the integration in D is performed then \ the integrations in
mutually complementary projections of D namely D$_{1}$ and D$_{2}$ are to be
performed along the lines parallel to x$_{1}$-axis in $\omega_{1}$ plane and
y$_{1}$-axis in $\omega_{2}$ plane respectively inside the strips $%
-\alpha<x_{2},y_{2}<\beta$ \ by using \ the contour integration technique.
So using $\ref{6.10}$ and $\ref{6.16},$ we obtain that

\begin{align*}
& \frac{1}{2\pi}\tint \limits_{D}\exp\{-i_{1}(\omega t)\widehat{f}%
(\omega)d\omega \\
& =\frac{1}{2\pi}\tint
\limits_{D}\exp\{-i_{1}(\omega_{1}e_{1}+\omega_{2}e_{2})t\}\{\widehat{f}%
_{1}(\omega _{1})e_{1}+\widehat{f}_{2}(\omega_{2})e_{2}\}d(\omega_{1}e_{1}+%
\omega_{2}e_{2}) \\
& =[\frac{1}{2\pi}\tint \limits_{D_{1}}\exp(-i_{1}\omega_{1}t)\widehat{f}%
_{1}(\omega_{1})d\omega_{1}]e_{1}+[\frac {1}{2\pi}\tint
\limits_{D_{2}}\exp(-i_{1}\omega_{2}t)\widehat{f}_{2}(\omega_{2})d%
\omega_{2}]e_{2} \\
& =[\frac{1}{2\pi}\tint
\limits_{-\infty+i_{1}x_{2}}^{\infty+i_{1}x_{2}}\exp(-i_{1}\omega_{1}t)%
\widehat{f_{1}}(\omega_{1})d\omega_{1}]e_{1}+[\frac {1}{2\pi}\tint
\limits_{-\infty+i_{1}y_{2}}^{\infty+i_{1}y_{2}}\exp(-i_{1}\omega_{2}t)%
\widehat{f_{2}}(\omega_{2})d\omega_{2}]e_{2} \\
& =f(t)e_{1}+f(t)e_{2}=f(t)\text{ }
\end{align*}
which guarantees the existence of Fourier inverse transform for
bicomplex-valued functions.

In the following, we define the bicomplex version of Fourier inverse
transform method.

\begin{definition}
Let $\widehat{f}(\omega)$ be the bicomplex Fourier transform of a real
valued continuous function $f(t)$ for $-\infty<t<\infty$ which is
holomorphic in $\ref{6.18}$.The Fourier inverse transform of $\widehat{f}%
(\omega)$\ is defined as%
\begin{equation*}
f(t)=\text{ }=\frac{1}{2\pi}\tint \limits_{D}\exp\{-i_{1}(\omega t)\widehat{f%
}(\omega)d\omega
\end{equation*}
where D is given by $\ref{6.20}$. On using $\ref{6.14}$,$\ref{6.15}$ and $%
\ref{6.17}$ this inversion method amounts to
\end{definition}

\begin{equation*}
f(t)=i_{1}e_{1}\tsum \limits_{Im(\omega_{2}^{(k)})>0}\func{Re}\text{s}%
\{\exp(-i_{1}\omega_{1}t)\widehat{f_{1}}(\omega
_{1}):\omega_{1}=\omega_{1}^{(k)}\}
\end{equation*}

\bigskip%
\begin{equation}
+i_{1}e_{2}\tsum \limits_{Im(\omega_{2}^{(k)})>0}\func{Re}\text{s}%
\{\exp(-i_{1}\omega_{2}t)\widehat{f_{2}}(\omega
_{2}):\omega_{2}=\omega_{2}^{(k)}\}\text{ for }t<0  \label{6.21}
\end{equation}

and

\begin{equation*}
f(t)=-i_{1}e_{1}\sum_{Im(\omega_{1}^{(k)})<0}\func{Re}\text{s}%
\{\exp(-i_{1}\omega_{1}t)\widehat{f_{1}}(\omega_{1}):\omega_{1}=\overline {%
\omega}_{1}^{(k)}\}
\end{equation*}

\begin{equation}
-i_{1}e_{2}\sum_{Im(\omega_{2}^{(k)})<0}\func{Re}\text{s}\{\exp
(-i_{1}\omega_{2}t)\widehat{f_{2}}(\omega_{2}):\omega_{2}=\omega_{2}^{(k)}\}%
\text{ for }t>0\text{.}  \label{6.22}
\end{equation}
We assign the value of $f(t)$ at $t=0$ fulfilling the requirement of
continuity of it in the whole real line $(-\infty<t<\infty).$

The following examples will make our notion clear:

\begin{example}
1.\textbf{\ }If \ \textbf{\ \ }$\widehat{f}(\omega)=\frac{2a}{a^{2}+\omega
^{2}}$\textbf{\ }for a\TEXTsymbol{>}0\ \textbf{\ }then%
\begin{equation*}
\widehat{f_{1}}(\omega_{1})=\frac{2a}{a^{2}+\omega_{1}^{2}},
\end{equation*}%
\begin{equation*}
\mathbf{\ \ }\widehat{f_{2}}(\omega_{2})=\frac{2a}{a^{2}+\omega_{2}^{2}}
\end{equation*}
\end{example}

and in each of $\omega _{1}$ and $\omega _{2}$ planes the poles are simple
at $\ i_{1}a$ and $i_{1}a$ . Now employing $\ref{6.21}$ and $\ref{6.22},$
for t \TEXTsymbol{<} 0 we\ obtain that 
\begin{equation*}
f(t)=i_{1}e_{1}\func{Re}\text{s}\{\exp (-i_{1}\omega _{1}t)\frac{2a}{%
a^{2}+\omega _{1}^{2}}:\omega _{1}=i_{1}a\}
\end{equation*}

\begin{equation*}
+i_{1}e_{2}\func{Re}\text{s}\{\exp(-i_{1}\omega_{2}t)\frac{2a}{%
a^{2}+\omega_{2}^{2}}:\omega_{2}=i_{1}a\}
\end{equation*}

\begin{equation*}
=i_{1}e_{1}\{-i_{1}\exp(at)\}+i_{1}e_{2}\{-i_{1}\exp(at)\}
\end{equation*}

\begin{equation*}
=\exp(-a|t|)
\end{equation*}

and for \ t \TEXTsymbol{>}0,%
\begin{equation*}
f(t)=-i_{1}e_{1}\func{Re}\text{s}\{\exp(-i_{1}\omega_{1}t)\frac {2a}{%
a^{2}+\omega_{1}^{2}}:\omega_{1}=i_{1}a\}
\end{equation*}

\begin{equation*}
-i_{1}e_{2}\func{Re}\text{s}\{\exp(-i_{1}\omega_{2}t)\frac{2a}{%
a^{2}+\omega_{2}^{2}}:\omega_{2}=i_{1}a\}
\end{equation*}

\begin{equation*}
=-i_{1}e_{1}\{i_{1}\exp(-at)\}-i_{1}e_{2}\{i_{1}\exp(at)\}
\end{equation*}

\begin{equation*}
=\exp(-a|t|).
\end{equation*}
Now for the continuity of t in \ the real line, we find $f(0)=1$. Thus the
Fourier inverse transform of $\widehat{f}(\omega)$ is $f(t)=\exp(-a|t|)$.

\begin{example}
\textbf{2. }If 
\begin{equation*}
\widehat{f}(\omega)=\frac{1}{2}[\frac{1}{\omega+\omega_{0}+\frac{i_{1}}{T}}-%
\frac{1}{\omega-\omega_{0}+\frac{i_{1}}{T}}]\text{ for }T,\omega_{0}>0
\end{equation*}
\end{example}

then in each of $\omega _{1}$ and $\omega _{2}$\ plane the poles are at ($%
\omega _{0}-\frac{i_{1}}{T})$ and ($-\omega _{0}-\frac{i_{1}}{T})$\ . For
both the poles the imaginary components are negative and so the poles are in
lower half of both the planes. In otherwords, no poles exist in upper half
of $\omega _{1}$ or $\omega _{2}$\ planes and as its consequence $f(t)=0$
for t \TEXTsymbol{<} 0. Now at t \TEXTsymbol{>} 0,%
\begin{equation*}
f(t)=-i_{1}e_{1}\func{Re}\text{s}\{\exp (-i_{1}\omega _{1}t)\widehat{f_{1}}%
(\omega _{1}):\omega _{1}=-\omega _{0}-\frac{i_{1}}{T}\}
\end{equation*}

\begin{equation*}
-i_{1}e_{1}\func{Re}\text{s}\{\exp(-i_{1}\omega_{1}t)\widehat{f_{1}}%
(\omega_{1}):\omega_{1}=\omega_{0}-\frac{i_{1}}{T}\}
\end{equation*}

\begin{equation*}
-i_{1}e_{2}\func{Re}\text{s}\{\exp(-i_{1}\omega_{2}t)\widehat{f_{2}}%
(\omega_{2}):\omega_{2}=-\omega_{0}-\frac{i_{1}}{T}\}
\end{equation*}

\begin{equation*}
-i_{1}e_{2}\func{Re}\text{s}\{\exp(-i_{1}\omega_{2}t)\widehat{f_{2}}%
(\omega_{2}):\omega_{2}=\omega_{0}-\frac{i_{1}}{T}\}
\end{equation*}

\begin{align*}
& =-i_{1}e_{1}\frac{1}{2}\exp (-\frac{t}{T})\exp (i_{1}\omega _{0}t) \\
& +i_{1}e_{1}\frac{1}{2}\exp (-\frac{t}{T})\exp (-i_{1}\omega _{0}t) \\
& -i_{1}e_{2}\frac{1}{2}\exp (-\frac{t}{T})\exp (i_{1}\omega _{0}t) \\
& +i_{1}e_{2}\frac{1}{2}\exp (-\frac{t}{T})\exp (-i_{1}\omega _{0}t) \\
& =-i_{1}\frac{1}{2}\exp (-\frac{t}{T})\exp (i_{1}\omega _{0}t) \\
& +i_{1}\frac{1}{2}\exp (-\frac{t}{T})\exp (-i_{1}\omega _{0}t) \\
& =\exp (-\frac{t}{T})\sin (\omega _{0}t).
\end{align*}%
Finally, the continuity of $f(t)$ in the whole real line implies that $%
f(0)=0 $.

\begin{equation*}
\text{- - - - - }\star\text{ - - - - -}
\end{equation*}

\newpage


\bigskip

\begin{thebibliography}{99}
\bibitem{Boc} S. Bochner and K. Chandrasekharan: Fourier transforms, Annals
of Mathematics Studies, Princeton University Press, Princeton, No. 19(1949).

\bibitem{Banerjee} A.Banerjee, S.K.Datta and Md. A.Hoque: Fourier transform
for functions of bicomplex variables, Asian journal of mathematic and its
applications,Vol.2015(2015),pp.1-18.

\bibitem{Gal} S.Gal: Introduction to geometric function theory of
hypercomplex variables, Nova Science Publishers, Vol.XVI (2002),pp.319-322.

\bibitem{Goyal} R.Goyal: Bicomplex polygamma function,Tokyo Journal of
Mathematics,Vol.30, No.2 (2007), pp.523-530

\bibitem{Ham} W. R. Hamilton: Lectures on quaternions containing a
systematic statement of a new mathematical method, Dublin, 1853.

\bibitem{Kai} G.Kaiser Birkhauser: A friendly guide to wavelets, Boston,1994.

\bibitem{Mot} A.Motter and M. Rosa : Hyperbolic calculus, Adv. Appl.
Clifford Algebra,Vol.8,No.1(1998),pp.109-128.

\bibitem{Mart} E.Martineau and D.Rochon : On a bicomplex distance estimation
for tetrabrot, Int. J. Bifurcation Chaos, Vol.15, No.9(2005),pp.3039-3035.

\bibitem{Mat} J.H.Mathews and R.W.Howell : Complex analysis for mathematics
and engineers, Narosa Publication , 2006.

\bibitem{Sid} Y.V.Sidorov, M.V.Fedoryuk and M.I.Shabunin : Lectures on the
theory of functions of complex variable, Mir Publishers,Moscow(1985).

\bibitem{Segre} C.Segre : Le rappresentazioni reali delle forme complesse a
gli enti iperalgebrici , Math.Ann.,Vol.40(1892), pp.413-467.
\end{thebibliography}
\end{document}